\documentclass[12pt]{article}
\usepackage{amsmath}
\usepackage{amsthm}
\usepackage{amsfonts}
\usepackage{amsfonts,amsmath,amstext,amsbsy,euscript,amssymb, graphics}
\usepackage{rotating}
\usepackage{enumerate}
\def\M{\mathcal{M}}
\def\P{\mathcal{P}}
\def\N{\mathcal{N}}
\def\B{\mathcal{B}}
\def\T{\mathcal{T}}
\def\W{\text{W}}
\def\RN{\text{RN}}
\def\RW{\text{RW}}

\def\Q{\cal{Q}}

\theoremstyle{plain}
\newtheorem{Thm}{Theorem}

\newtheorem{Lemma}[Thm]{Lemma}

\newtheorem{Cor}[Thm]{Corollary}

\newtheorem{Ex}[Thm]{Example}
\newtheorem{Rem}[Thm]{Remark}
\newtheorem{thm}{Theorem}
\begin{document}
\title{Rational Heap Games}
\author{Urban Larsson,\\ Mathematical Sciences, \\
  Chalmers University of Technology and University of Gothenburg}
\maketitle
\begin{abstract}
We study variations of classical combinatorial games on two finite heaps of tokens, a.k.a. \emph{subtraction games}. Given non-negative integers $p_1,q_1, p_2,q_2$, where $p_1q_2 > q_1p_2$, $p_1>0$ and $q_2>0$, two players alternate in removing $(m_1,m_2)\ne (0,0)$ tokens from the respective heaps, where the allowed ordered pairs of non-negative integers are given by a certain move set $(m_1,m_2)\in\M$. There is a restriction imposed on the allowed heap sizes $(X, Y)$, they must satisfy $Xq_1\le Yp_1$ and $Yp_2\le Xq_2$.  A player who cannot move loses and the other player wins. For a certain restriction of these games, namely where each allowed move option $(m_1,m_2)$ is of the form $(sp_1+tp_2,sq_1+tq_2)$, for some ordered pair of non-negative integers $(s,t)\ne (0,0)$, we show that all games have equivalent outcomes via a certain surjective map to a canonical subtraction game. Other interests in our games are various interactions with classical combinatorial games such as \emph{Nim {\rm and} Wythoff Nim}.
\end{abstract}

\section{Introduction}\label{S:1}
We study \emph{impartial} \emph{subtraction} games on two \emph{heaps of tokens} \cite{G1966, DR2010, L1, L2, LHF2011, L2012, LW}. For a background, see also \cite{BCG1982}. There are two players who obey the same rules and alternate in moving. We follow the \emph{normal play} convention, meaning that a player who cannot move loses and the other player wins. 

In general an impartial game consists of a \emph{position}, which contains information of whose turn it is and describes the given state of the game, and a \emph{ruleset}, which decides what move options there are for a given position. Sometimes the term ``game'' is adapted to mean only the ruleset, other times the term ``position'' encompasses the usual meaning of a game, that is whenever the ruleset is understood. Here we study so-called \emph{invariant} games \cite{G1966, DR2010, L1, L2, LHF2011, L2012, LW}, where the rules in general do not depend on the given position. For example, in \emph{2-heap Nim} \cite{B1902} the players remove tokens from precisely one of two finite heaps, at least one token and at most a whole heap. The game is invariant in the sense that, independent of the size of a heap, a given number of tokens can be removed provided the heap contains at least this number of tokens. As for all impartial games, this game has a \emph{perfect winning strategy}, here: the second player to move wins if and only if the two heaps have the same non-negative number of tokens. 

In general, for an impartial game $G$ without drawn moves, the \emph{outcome}, of a given position, is either a \emph{previous player} win ($\P$) or a \emph{next player} win ($\N$). It belongs to $\P=\P(G)$, if none of its options is in $\P$. Otherwise it belongs to $\N$. This gives a recursive characterization of all $\P$-positions of a given game $G$, starting with the \emph{terminal} positions in $\T(G)\subseteq \P$, from which no move is possible. 

\subsection{Rational heap games}
Let 
\[ Q:=\left(
\begin{array}{cc}
    p_1 & q_1 \\
    p_2 & q_2 \\
\end{array} \right),\]
where $p_1>0, q_1, p_2, q_2>0$ are given non-negative integer \emph{game constants}, with $\det Q = p_1q_2 - q_1p_2 > 0$. We let the allowed \emph{positions} or \emph{heap-sizes} be represented by ordered pairs of non-negative integers, bounded by `rational slopes' as given by the set 
\begin{align}\label{2}
\B_Q:=\left\{(X,Y)\mid Xq_1\le Yp_1\text{ and }Yp_2\le Xq_2\right\}. 
\end{align} 
In particular, for the special cases $q_1=p_2=0$ all combinations of heap-sizes are allowed; we omit the index and denoted this set simply by 
\begin{align}\label{B}
\B:=\left\{(X,Y)\mid 0\le X, 0\le Y\right\}.
\end{align} 
Also, we let $\B':=\B\setminus \{(0,0)\}$.
Following \cite{DR2010,LHF2011,L2012}, a two heap \emph{subtraction game}, $G=G(\M)$, is defined via a given set of ordered pairs of non-negative integers $\M \subseteq \B'$. A legal move from $(X,Y)\in \B$ is to $(X-m_1, Y-m_2)$ for some $(m_1, m_2)\in \M$ and provided $(X-m_1,Y-m_2)\in \B$.

Let $(X,Y)\in \B_Q$. In general, each allowed move for a \emph{$\B_{\Q}$-subtraction game} is of the form
\begin{align}\label{3}
(X,Y)\rightarrow (X-m_1, Y-m_2), 
\end{align}
where $(m_1,m_2)$ belongs to a given set of ordered pairs of non-negative integers $\M\subseteq\B'$, and provided $(X-m_1,Y-m_2)\in \B_{\Q}$. Our main interest in this paper is the following subset of the $\B_{\Q}$-subtraction games. For a purpose that will become clear later (in Lemma \ref{L:5}) we will alter the notation somewhat.

In a \emph{$\Q$-subtraction game} $G_Q=G_Q(\M)$ we require that $(m_1,m_2)$ in (\ref{3}) is of the form $m_1 = p_1s + p_2t$ and $m_2 = q_1s + q_2t$ where $(s,t)\in \M\subseteq \B'$.  Hence for these games it is $(s,t)$ which belongs to ``the set of allowed moves'' $\M$ (and not $(m_1,m_2)$). Thus, a typical move in $G_Q(\M)$ is 
\begin{align}\label{4}
(X,Y)\rightarrow (X-p_1s-p_2t, Y-q_1s-q_2t), 
\end{align}
where $(s,t)\in \M$ and $(X-p_1s-p_2t, Y-q_1s-q_2t)\in \B_Q$. As before, whenever a $\Q$-subtraction game is a subtraction game we write simply $G(\M)$, that is whenever 
\[ Q=\left(
\begin{array}{cc}
    1 & 0 \\
    0 & 1 \\
\end{array} \right)\]
and the allowed heap sizes are as given by $\B$ in (\ref{B}).

Let us define a surjective map $\varphi_Q$, which takes as input a position in $\B_Q$ and produces as output a position in $\B$, 
\begin{align}\label{varphi}
\varphi_Q(X,Y) = \left(\left\lfloor \frac{Xq_2-Yp_2}{\det Q}\right\rfloor,\left\lfloor\frac{Yp_1-Xq_1}{\det Q}\right\rfloor\right),
\end{align} 
by $q_1X/p_1\le Y\le q_2X/p_2$. Our main result is that any $\Q$-subtraction game $G_Q:=G_Q(\M)$ is ``$\varphi_Q$-equivalent'' to the subtraction game $G:=G(\M)$ in the following sense. 

\begin{Thm}\label{T:1}
Given game constants $p_i, q_i$, suppose that $(X,Y)\in \B_Q$. Then 
$(X,Y)\in \P(G_Q)$ if and only if $\varphi_Q(X,Y)\in \P(G)$.
This is equivalent to $(A,B)\in \P(G)$ if and only if, for all $(x,y)\in \T_Q$, $(x+Ap_1+Bp_2, y+Aq_1+Bq_2)\in \P(G_Q)$. 
\end{Thm}

By this result we make the following definition: the $\Q$-subtraction games $G_{Q}(\M)$ and $G_{R}(\mathcal{L})$ are \emph{$\varphi$-equivalent} if $\M=\mathcal{L}$. The \emph{$\M$-canonical game} is the subtraction game $G(\M)$.

In Section \ref{S:2} we prove Theorem \ref{T:1}. A reader who wishes to study some examples before plunging into the proof of the main theorem should skip to Section \ref{S:3}, where we illustrate Theorem \ref{T:1} via generalizations of Nim and Wythoff Nim. This discussion is continued in Section \ref{S:4} with some open questions where we relate a certain ``splitting behavior'' of Wythoff type $Q$-subtraction games to similar $\B_{\Q}$-subtraction games.

\section{Proof of Theorem \ref{T:1}}\label{S:2}
A generic game has several terminal positions from which no move is possible.

\begin{Lemma}\label{L:1}
Given $\M\subseteq \B'$, the set
\begin{align*}
\T_Q := \left\{(x,y)\mid p_1(y-q_2)<q_1(x-p_2) \text{ and } p_2(y-q_1)>q_2(x-p_1)\right\}\subseteq \B_Q
\end{align*} 
is a subset of all \emph{terminal positions} of the $\Q$-subtraction game $G_Q(\M)$. In particular $\T_Q = \T(G_Q)$ is the set of all terminal positions, if and only if $\{(0,1),(1,0)\}\subseteq \M$. Also $\det Q=|\T_Q|$.
\end{Lemma}
\noindent{\bf Proof.} Since, by assumption $(x,y)\in \T_Q\subset \B_Q$, we need to show that $(x-sp_1-tp_2, y-sq_1-tq_2)\not\in \B_Q$ for all positive integers $s$ and $t$. Hence, by definition of $\B_Q$, we need to show that $(x-sp_1-tp_2)q_1> (y-sq_1-tq_2)p_1$ and that $(x-sp_1-tp_2)q_2 < (y-sq_1-tq_2)p_2$ for all positive integers $s$ and $t$, but this is clear since, by definition of $\T_Q$, it holds for $s=1$ and $t=1$, and the game constants satisfy $p_1q_2 > q_1p_2$ . 

The set $\T_Q$ is the complete set of terminal positions of $G_Q(\M)$ if $$\{(0,1),(1,0)\}~\subseteq~\M$$ since $(x,y)\not\in \T_Q$ gives $p_1(y-q_2)\ge q_1(x-p_2)$ or $p_2(y-q_1)\le q_2(x-p_1)$. Hence  $(x-p_2, y-q_2)\in \B_Q$ or $(x-p_1, y-q_1)\in \B_Q$ respectively. For ``only if'', suppose that $(0,1)\not\in \M$, then $(x+p_2,y+q_2)\in \T(G_Q)\setminus \T_Q$. The other case is similar. 
\hfill $\Box$\\

The next two lemmas discuss how the positions in $\B_Q$ can be viewed as linear translations of those in $\T_Q$.

\begin{Lemma}\label{L:2}
Let $(x,y)\in \T_Q$. Then $(x+Ap_1+Bp_2,y+Aq_1+Bq_2)\in \B_Q$ if and only if $(A,B)\in \B$.
\end{Lemma}

\noindent{\bf Proof.} What is required is to show that $$(x+Ap_1+Bp_2)q_1\le (y+Aq_1+Bq_2)p_1$$ and 
\begin{align}\label{second}
(x+Ap_1+Bp_2)q_2\ge (y+Aq_1+Bq_2)p_2 
\end{align}
if and only if both $A$ and $B$ are non-negative, assuming that $(x,y)\in \T_Q\subset \B_Q$. But, as in the proof of Lemma \ref{L:1}, the ``if''-part is immediate by $p_1q_2 > q_1p_2$.

Hence, suppose that $A<0$.  By definition there is no move from any $\Q$-subtraction game from $(x,y)\in \T_Q$. Hence, by negativity of $A$, we must have $(x+Ap_1)q_2 < (y+Aq_1)p_2$, which contradicts (\ref{second}). The case $B<0$ is similar. \hfill $\Box$\\

\begin{Lemma}\label{L:3}
If $(X,Y)\in \B_Q$, then $(X,Y) = (x+Ap_1+Bp_2, y+Aq_1+Bq_2)$, for some unique $(x, y)\in \T_Q$ and some unique non-negative integers $A$ and $B$. 
\end{Lemma}

\noindent{\bf Proof.} Suppose that $(X,Y) = (x+Ap_1+Bp_2, y+Aq_1+Bq_2)=(x'+A'p_1+B'p_2, y'+A'q_1+B'q_2)$, with $(x',y')\in \T_Q$, $x\ge x', y\ge y'$ and non-negative integers $A'$ and $B'$. Suppose that $x>x'$. This gives $x-x'=(A'-A)p_1+(B'-B)p_2>0$ which implies $y-y' = (A'-A)q_1+(B'-B)q_2 > 0$ which is impossible by Lemma \ref{L:1}. Hence $(A'-A)p_1=(B-B')p_2$ and $(A'-A)q_1=(B-B')q_2$ which gives $q_1p_2=p_1q_2$ which is impossible.

Next, let us find $x,y,A,B$ such that $(X, Y) = (x+Ap_1+Bp_2, y+Aq_1+Bq_2)$. We have that $Xq_1 \le Yp_1$ and $Yp_2 \le Xq_2$. Hence, by $q_1p_2 < p_1q_2$, there is a largest non-negative $B$ such that $(X - Bp_2)q_1 \le (Y - Bq_2)p_1$ and a largest non-negative $A$ such that $(Y - Aq_1)p_2 \le (X - Ap_1)q_2$. Then $(X-Ap_1-Bp_2,Y-Aq_1-Bq_2)\in \T_Q$ defines $(x,y)\in \T_Q$ by our choices of $A$ and $B$.
\hfill $\Box$\\

\begin{Rem}\label{R:1} 
Following Lemma \ref{L:3}, we say that (the game) $(X,Y)=(x+Ap_1+Bp_2, y+Aq_1+Bq_2)$ belongs to the \emph{$(x,y)$-class}, where $(x,y)\in \T_Q$. Then, by $(X,Y) = (A,B)Q$, $(X,Y)$ belongs to the $(0,0)$-class if and only if the associated restriction of $\varphi_Q$ is $\varphi_Q|_{(0,0)}(X,Y)=(X,Y)Q^{-1} = (A,B)$.
\end{Rem}

We get the following consequence of the above lemmas for our $Q$-subtraction games.

\begin{Lemma}\label{L:5}
Let $\M\subseteq \B'$. Suppose that there is a move  of the form $(A,B)\rightarrow (A-s,B-t)$ in the subtraction game $G(\M)$, then for any given game constants $p_i,q_i$ and for each $(x,y)\in \T_Q$, there is a move in the $Q$-subtraction game $G_Q(\M)$ of the form 
\begin{align}\label{moveQ}
&(x+Ap_1+Bp_2,y+Aq_1+Bq_2)\rightarrow \notag\\&(x+(A-s)p_1+(B-t)p_2,y+(A-s)q_1+(B-t)q_2). 
\end{align}
Suppose on the other hand that there is a move in $G_Q$ from $(X,Y)\in \B_Q$ via $(s,t)\in \M$. Then there is a corresponding move $(A,B)\rightarrow (A-s,B-t)$ in $G$ where $(X,Y)=(x+Ap_1+Bp_2,y+Aq_1+Bq_2)$ for some $(x,y)\in \T_Q$.
\end{Lemma}

Thus, with notation as in Lemma \ref{L:5} and Remark \ref{R:1}, each game in the $(x,y)$-class ends in $(x,y)$ if $\{(0,1),(1,0)\}\subseteq \M$.

Let us restate and prove our main theorem, where the function $\varphi_Q$ is as in (\ref{varphi}).

\begin{thm}
Given game constants $p_i, q_i$, suppose that $(X,Y)\in \B_Q$. Then 
$(X,Y)\in \P(G_Q)$ if and only if $\varphi_Q(X,Y)\in \P(G)$.
This is equivalent to $(A,B)\in \P(G)$ if and only if, for all $(x,y)\in \T_Q$, $(x+Ap_1+Bp_2, y+Aq_1+Bq_2)\in \P(G_Q)$. 
\end{thm}

\noindent{\bf Proof.} We begin with the second part. Suppose that $(A,B)\in \P(G)$. Then, none of its options is in $\P(G)$. We need to prove that none of the options in the $Q$-subtraction game $G_Q$ from $(x+Ap_1+Bp_2, y+Aq_1+Bq_2)$ is in $\P(G_Q)$. We have that, for all $(s,t)\in \M$ such that $A-s\ge 0$ and $B-t\ge 0$, $(A-s,B-t)\in \N(G)$. Then, by Lemma \ref{L:5}, induction gives that $(x+(A-s)p_1+(B-t)p_2, y+(A-s)q_1+(B-t)q_2)\in \N(G_Q)$ if and only if $A-s\ge 0$ and $B-t\ge 0$.

If, on the other hand $(A, B)\in \N(G)$, then there is an option $(A-s, B-t)\in \P(G)$, with $(s, t)\in \M$, and so induction gives that $(x+(A-s)p_1+(B-t)p_2, y+(A-s)q_1+(B-t)q_2)\in \P(G_Q)$.

For the first part, by Lemma \ref{L:3}, if $(X,Y)\in \B_Q$, then $(X,Y) = (x+Ap_1+Bp_2, y+Aq_1+Bq_2)$, for some unique $(x, y)\in \T_Q$ and some unique non-negative integers $A$ and $B$. We plug this into the expression

 \begin{align*}
\varphi_Q(X,Y) &= \left(\left\lfloor \frac{Xq_2-Yp_2}{\det Q}\right\rfloor,\left\lfloor\frac{Yp_1-Xq_1}{\det Q}\right\rfloor\right)\\
&=\left(\left\lfloor \frac{xq_2-yp_2+A\det Q}{\det Q}\right\rfloor,\left\lfloor\frac{yp_1-xq_1+B\det Q}{\det Q}\right\rfloor\right)
&=(A, B),
\end{align*}
by $0\le \frac{yp_1-xq_1}{\det Q}<1$ and $0\le \frac{xq_2-yp_2}{\det Q}<1$, since $(x,y)\in \T_Q$. But then the first part of the proof gives the result.
\hfill $\Box$\\

\section{Examples}\label{S:3}
Our first example of a 2-heap subtraction game, 2-heap Nim, was discussed briefly in the introduction.

The game of \emph{Wythoff Nim}  \cite{W1907} is also played on two heaps, all the moves in Nim are allowed and also the possibility of removing the same positive number of tokens from both heaps in one and the same move, a number bounded by the number of tokens in the smallest heap. 

\begin{figure}[ht!]
\begin{center}
\includegraphics[width=0.4\textwidth]{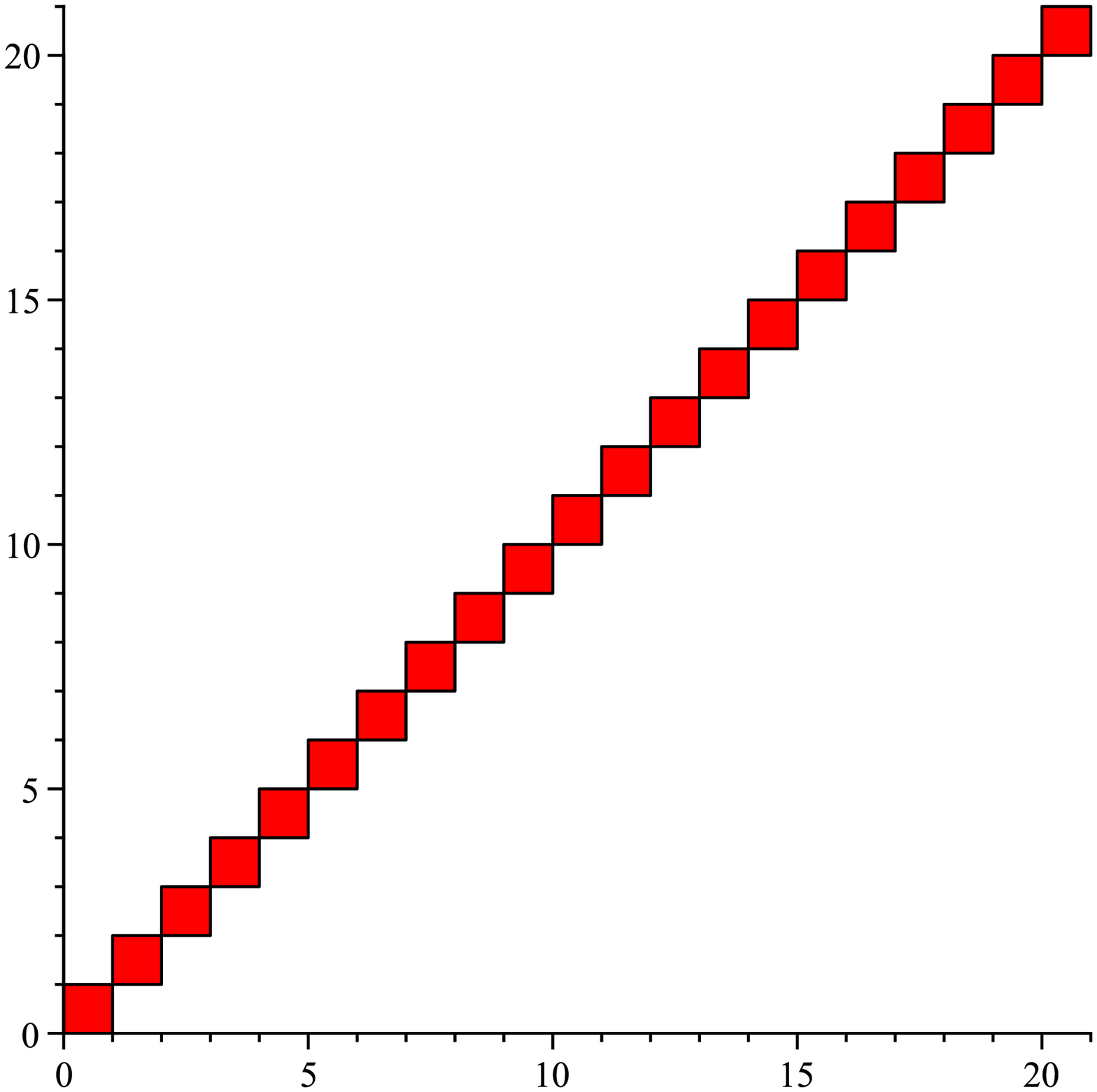}\hspace{0.5 cm}
\includegraphics[width=0.4\textwidth]{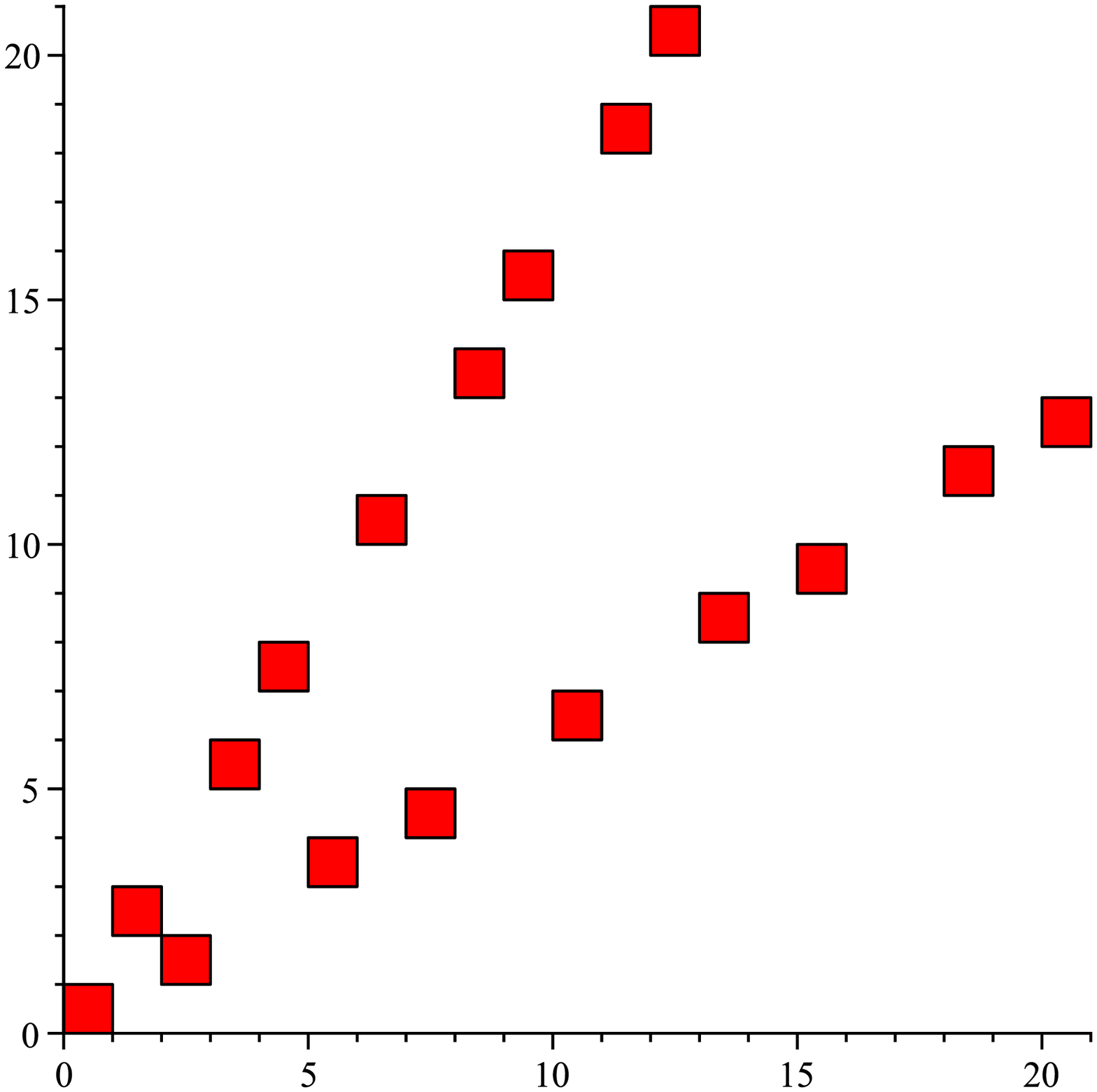}
\end{center}\caption{The red squares represent the initial patterns of $\P$-positions of 2-heap Nim, $(0,0),(1,1),\ldots$ to the left and Wythoff Nim, $(0,0),(1,2),(2,1),(3,5),(5,3),\ldots$ to the right.}\label{F:Nim} 
\end{figure}

See Figure \ref{F:Nim} 
for the $\P$-positions of Nim and Wythoff Nim (\W) respectively. For the latter game, it is known \cite{W1907} that they are of the form: 

\begin{align}\label{phi}
\P(\W)=\{(X,Y),(Y,X)\mid X=\lfloor \phi n\rfloor , Y=\lfloor \phi^2 n \rfloor\},
\end{align} 
where $n$ ranges over the non-negative integers and $\phi := \frac{1+\sqrt{5}}{2}$ denotes the golden ratio. 

\begin{figure}[ht!]
\begin{center}
\includegraphics[width=0.5\textwidth]{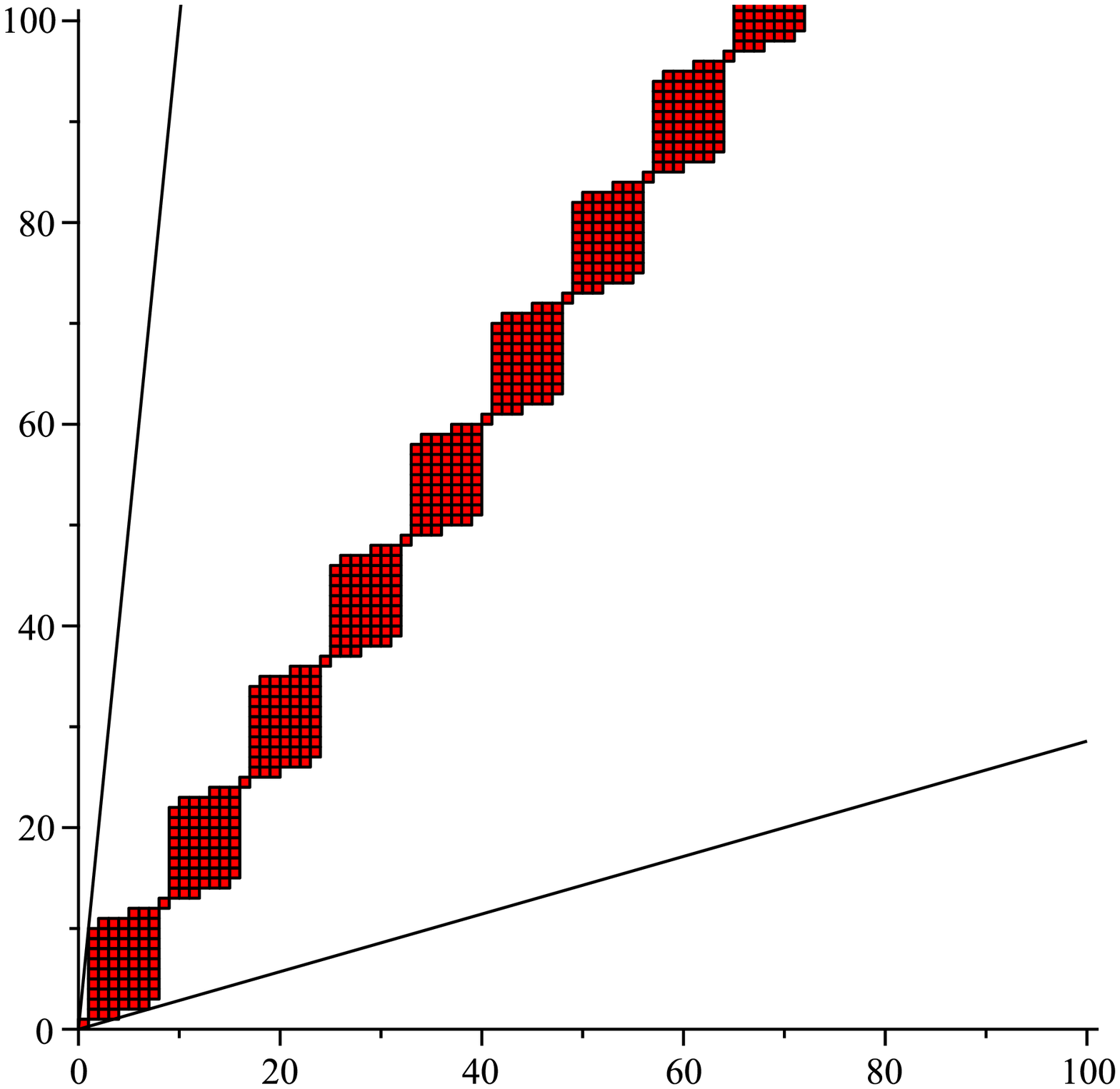}
\end{center}\caption{The red squares represent the initial $\P$-positions for the game $(\frac{2}{7},\frac{10}{1})$-RN. The two lines represent the bounds for the the ratio of the heap-sizes.}\label{F:1} 
\begin{center}
\includegraphics[width=0.5\textwidth]{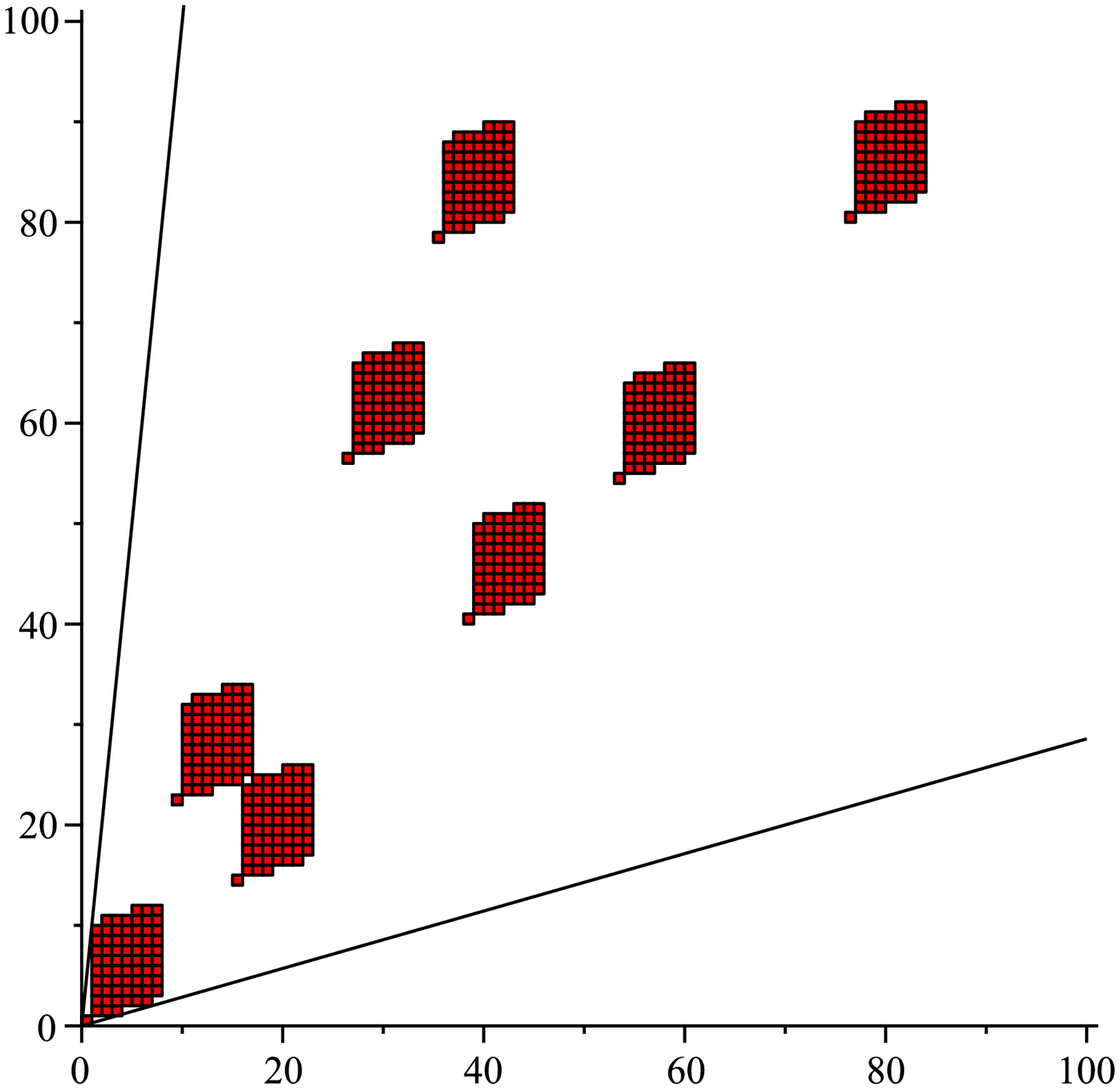}
\end{center}\caption{The red squares represent all initial $\P$-positions for the game $(\frac{2}{7},\frac{10}{1})$-RW.}\label{F:2} 
\end{figure}

\begin{Ex}
For our first example of a $Q$-subtraction game, a player may move 
\begin{align}\label{5}
(x,y)\rightarrow (x-p_it, y-q_it), 
\end{align}
for $i\in \{1,2\}$ and any positive integer $t$, provided $(x-p_it,y-q_it)\in \B_Q$. 

We denote this game by \emph{Rational Nim} (RN), $(\frac{q_1}{p_1},\frac{q_2}{p_2})$-RN. Here the ratios in the prefix are merely symbols and so for example $(\frac{1}{2},\frac{2}{3})$-RN and $(\frac{2}{4},\frac{4}{6})$-RN denote different games. Note that, given $Q$, this game may be represented simply by the set $\M = \{(0,t),(t,0)\mid t>0\}$ and that the case $p_1=q_2=1$, $p_2=q_1=0$ is the classical game of Nim on 2 heaps. See Figure \ref{F:1} for the $\P$-positions of an RN game.
\end{Ex}

An \emph{extension} of a game $G$ has the same set of positions as $G$, contains all the moves in $G$ and possibly some new. Thus a trivial extension of 2-heap Nim is Nim itself. A non-trivial extension is for example Wythoff Nim. Further, a \emph{$Q$-extension} is an extension of a $Q$-subtraction game which is also a $Q$-subtraction game. Thus another way of expressing this is that $G_Q(\M'')=G_Q(\M\cup\M')$ is a $Q$-extension of both $G_Q(\M)$ and $G_Q(\M')$. Note that $Q$ is fixed.

\begin{Ex}
Given game constants $p_i, q_i$, we denote by \emph{Rational Wythoff Nim} (RW), \emph{$(\frac{q_1}{p_1},\frac{q_2}{p_2})$-\RW}, the following $Q$-extension of $(\frac{q_1}{p_1},\frac{q_2}{p_2})$-RN. In addition to the moves in (\ref{2}), the new moves are of the form $$(x,y)\rightarrow (x-t(p_1+p_2),y-t(q_1+q_2))$$ provided that $(x-t(p_1+p_2),y-t(q_1+q_2))\in \B_Q$. Hence, given $Q$, this game may equivalently be represented by $\M=\{(0,t),(t,0),(t,t)\mid t>0\}$ and the classical game of Wythoff Nim is the game where $p_1=q_2=1$, $p_2=q_1=0$.
See Figure \ref{F:2} for the initial $\P$-positions of an RW game. 
\end{Ex}

Theorem \ref{T:1} gives that $\P(\RN)$ is periodic with period $(p_1+p_2,q_1+q_2)$ and $\P(\RW)$ is aperiodic similar to $\P(\W)$, as described by the following result. 
\clearpage
\begin{Cor}\label{C:1}
The $\P$-positions of Rational Nim are given by $\P(RN)=\{(x+n(p_1+p_2),y+n(q_1+q_2)\}$, where $(x,y)\in \T_Q$ and where $n$ ranges over the non-negative integers.

The set $\P(\RW)$ is given by all positions of the forms
\begin{align}
(x+p_1\lfloor \phi^2 n\rfloor+p_2\lfloor \phi n\rfloor,y+q_1\lfloor \phi^2 n\rfloor+q_2\lfloor \phi n\rfloor)
\end{align}
and 
\begin{align}
(x+p_1\lfloor \phi n\rfloor+p_2\lfloor \phi^2 n\rfloor,y+q_1\lfloor \phi n\rfloor+q_2\lfloor \phi^2 n\rfloor),
\end{align}
where $(x,y)\in \T_Q$ and $n$ is a non-negative integer.
\end{Cor}
 See also Corollary \ref{C:2} in the next section.

\section{When a game extension splits a set of $\P$-positions}\label{S:4}
Is it true that adjoining new moves to an existing impartial game changes its set of $\P$-positions? For the particular games studied in Section \ref{S:3} this is certainly true by Corollary \ref{C:1}, as is also illustrated in Figures $\ref{F:1}$ and $\ref{F:2}$. In fact, we have seen a particular shift of behavior inherited from the relation between Nim and Wythoff Nim. In describing RN's and RW's asymptotic behavior we go from one accumulation point to two accumulation points. Following the terminology in \cite{L1}, we note that the particular $Q$-extension of RN we have introduced, RW, \emph{splits} the `old' set of $\P$-positions of RN into two `new' $\P$ sets for which the ratios of the heap-sizes \emph{converge} to two distinct real numbers. We have the following simple consequence of Corollary~\ref{C:1}.

\begin{Cor}\label{C:2}
Let the game constants $p_i,q_i$ be given. 
Let $\P(\RN) = \{(A_n,B_n)\}$, with the $(A_n,B_n)$s in lexicographical order. Then $$\frac{B_n}{A_n}\rightarrow \frac{q_1+q_2}{p_1+p_2},$$ as $n\rightarrow \infty$. 
Let $\P(\RW) = \{(A_n^l, B_n^l),(A_n^u,B_n^u)\}$, where the \emph{lower} $\P$-positions are the $(A_n^l,B_n^l)$s in lexicographical order and the \emph{upper} $\P$-positions are the $(A_n^u,B_n^u)$s in lexicographical order for which $$\frac{B_n^l}{A_n^l}< \frac{q_1+q_2}{p_1+p_2}\le \frac{B_n^u}{A_n^u}.$$ Then 
$$\frac{B_n^l}{A_n^l}\rightarrow \frac{q_1+\phi(q_1+q_2)}{p_1+\phi(p_1+p_2)},$$ 
$$\frac{B_n^u}{A_n^u}\rightarrow \frac{q_2+\phi(q_1+q_2)}{p_2+\phi(p_1+p_2)},$$
as $n\rightarrow\infty$

\end{Cor}

\begin{figure}[ht!]
\begin{center}
\includegraphics[width=0.45\textwidth]{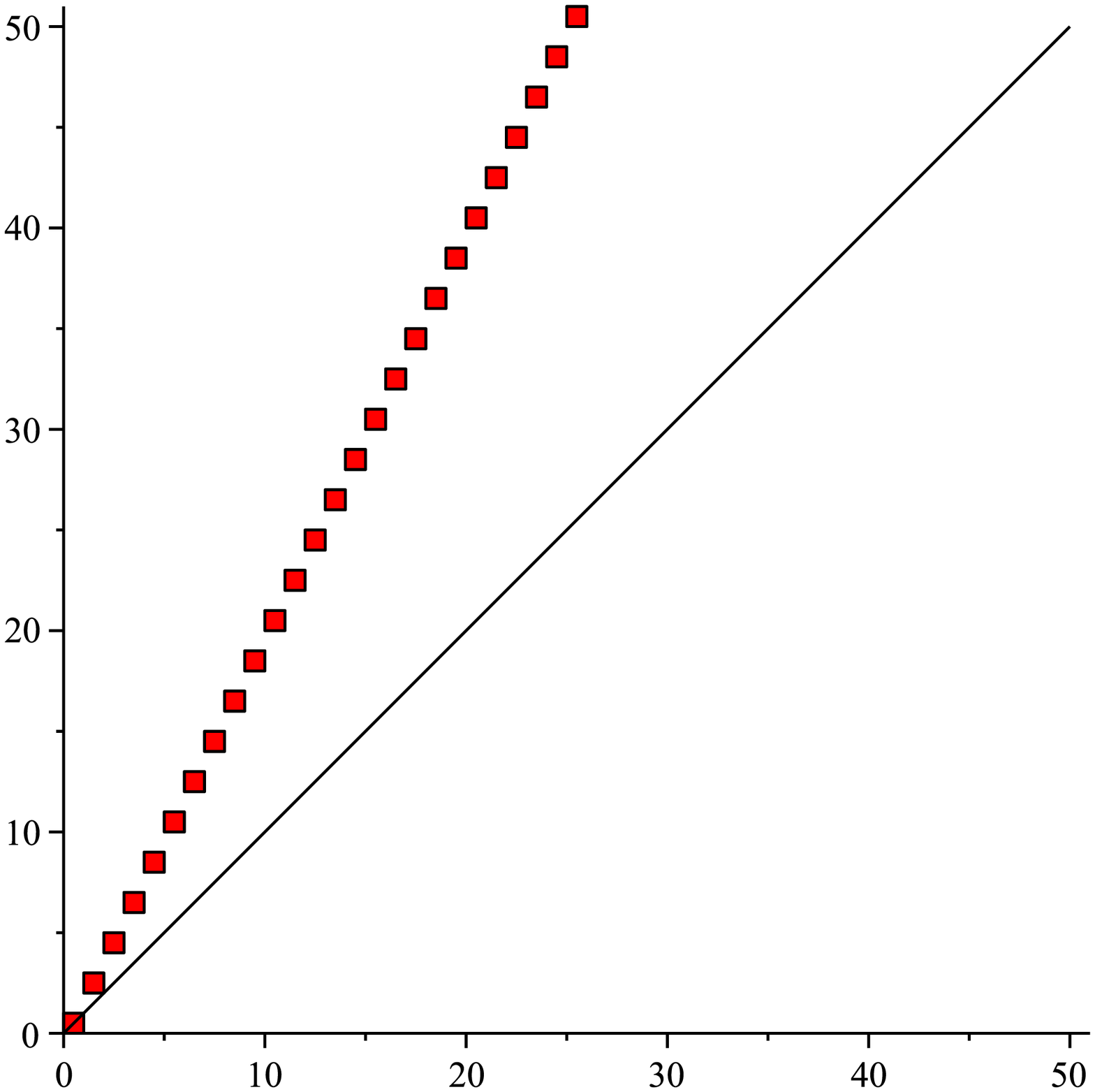}\hspace{0.3 cm}
\includegraphics[width=0.45\textwidth]{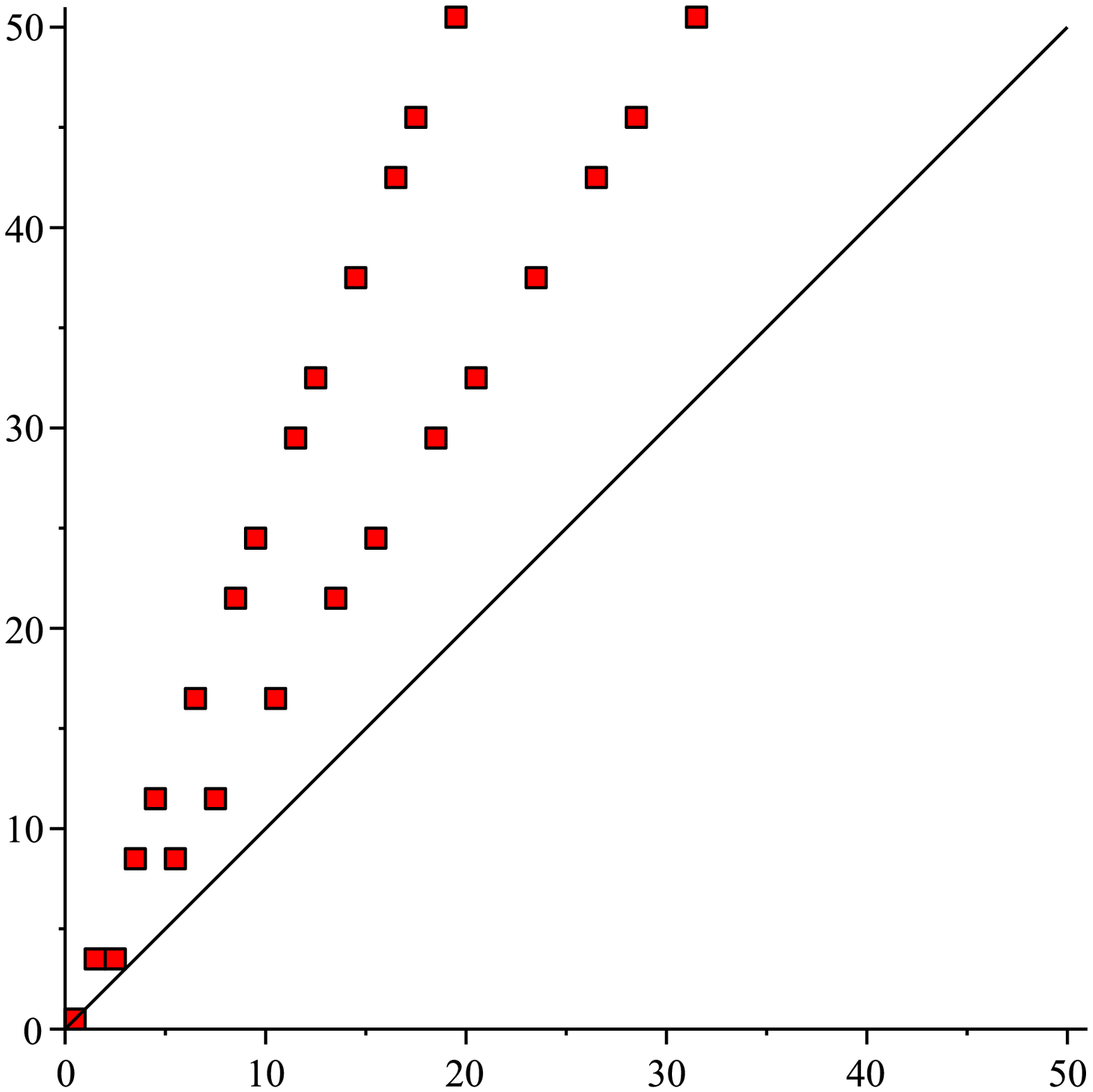}
\end{center}\caption{The red squares represent the initial $\P$-positions for the games $(\frac{1}{1},\frac{1}{0})$-RN and $(\frac{1}{1},\frac{1}{0})$-RW respectively.}\label{F3} 
\end{figure}

For another example, take the set $\mathcal{P}((\frac{1}{1},\frac{1}{0})\text{-\RW})$. The first few $\P$-positions are: $(0,0),(1,3),(2,3),(3,8),(5,8),(4,11),(7,11),\ldots .$  For all non-negative integers $n$, by Corollary \ref{C:1} and since the only terminal position is $(0,0)$, $$\P\left(\left(\frac{1}{1},\frac{1}{0}\right)\!-\!\RW\right) = \left\{(\lfloor \phi n\rfloor, \lfloor \phi n\rfloor+\lfloor \phi^2 n \rfloor), (\lfloor \phi^2 n \rfloor, \lfloor \phi n\rfloor+\lfloor \phi^2 n \rfloor)\right\}$$ and thus, the two convergents of the $\P$-positions of this game are $\phi$ and $\phi^2$, whereas for $(\frac{1}{1},\frac{1}{0})$-RN the single convergent is $2$. See also Figure \ref{F3}. This game has a particular interest as being a new simple `restriction' of the game $(1,2)$GDWN in \cite{L1}, the latter which is an extension of Wythoff Nim where the new moves are of the form: in one and the same move remove $t$ tokens from one of the heaps and $2t$ from the other, at most a whole heap. The game $(1,2)$GDWN is conjectured to split the upper (and lower) $\P$-positions of Wythoff Nim from the single convergent $\phi$ to a pair of convergents $1.478\ldots$ and $2.247\ldots$. See also \cite{L2011} for another example where a further split of the $\P$-positions of Wythoff Nim is obtained via a certain \emph{blocking maneuver} on the regular Wythoff Nim moves. 

Let us finish off this paper with some extensions of our sample game in Figure \ref{F:1}, that (unlike the game in Figure \ref{F:2}) are not $Q$-extensions, but indeed $\B_Q$-subtraction games. 
Our final three extensions of $(\frac{2}{7},\frac{10}{1})$-RN are as follows: adjoin moves of the form \\

\noindent game (a): $(X,Y)\rightarrow (X-4t,Y-6t)$,\\ 

\noindent game (b): $(X,Y)\rightarrow (X-4t,Y-4t)$,\\

\noindent game (c): $(X,Y)\rightarrow (X-8t,Y-4t)$,\\ 

for positive integers $t$, bearing in mind that the positions satisfy (\ref{2}). The first two extensions have complicated patterns of $\P$-positions which we do not yet understand, of which at least the first appears to exhibit a similar splitting `behavior' as does $(\frac{2}{7},\frac{10}{1})$-RW, but game (c) has the same set of $\P$-positions as does $(\frac{2}{7},\frac{10}{1})$-RN, which can be proved by elementary methods, see also \cite{FL1991, DFNR2010, L1, L2011} for related results. Hence the answer to the first question in this final section is negative; see also \cite{LW} for a discussion of this question in the context of heap games, computational complexity and algorithmic undecidability. 

\section{Discussion}
\begin{figure}[ht!]
\begin{center}
\includegraphics[width=0.5\textwidth]{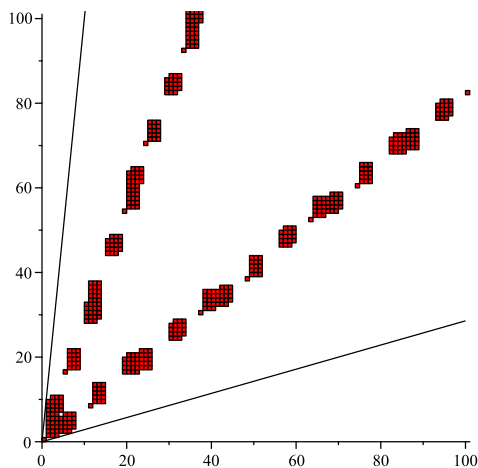}
\end{center}\caption{The red squares represent the initial $\P$-positions for extension (a) of $(\frac{2}{7},\frac{10}{1})$-RN.}\label{F:4} 
\end{figure}

\begin{figure}[ht!]
\begin{center}
\includegraphics[width=0.5\textwidth]{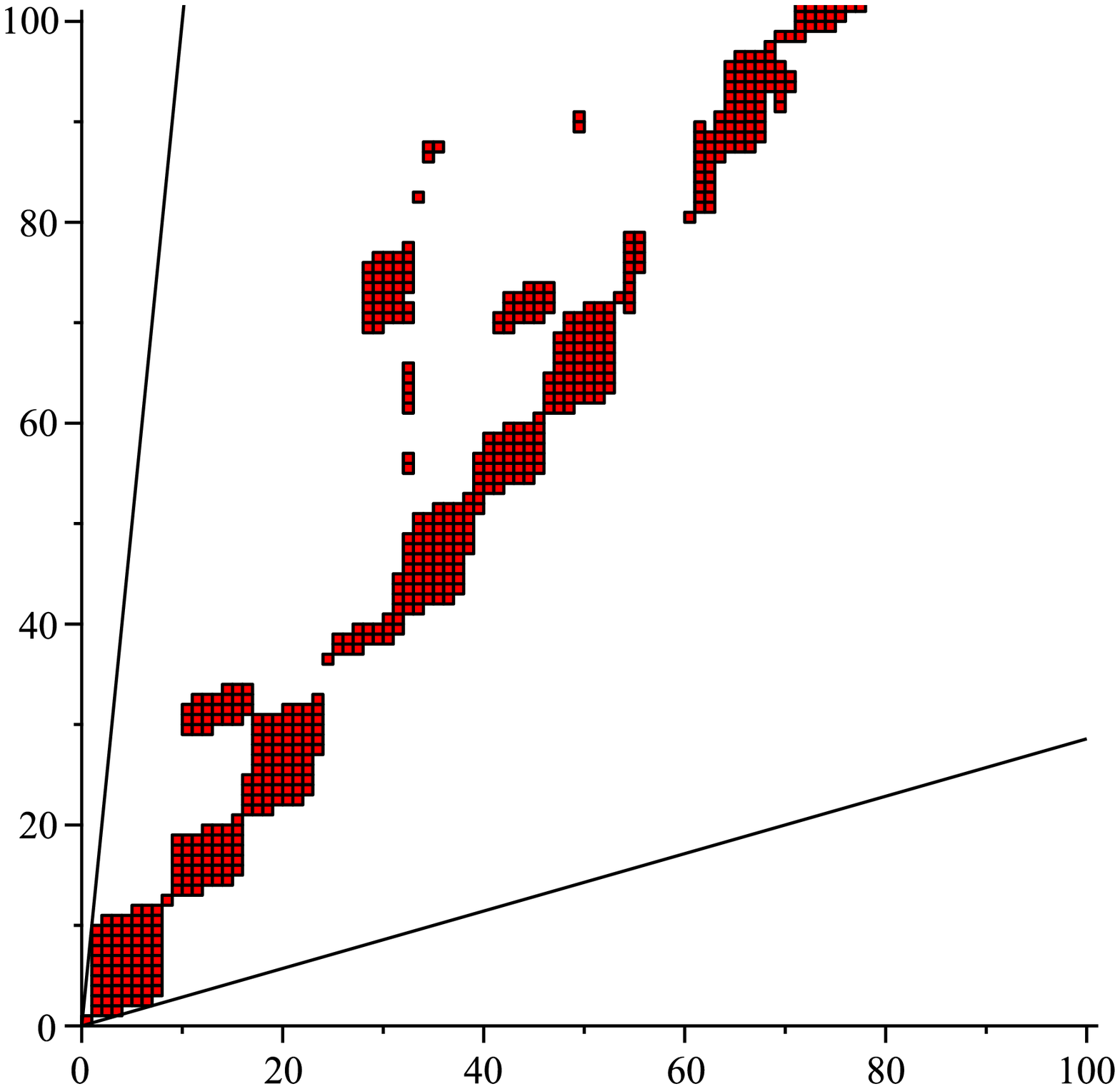}
\end{center}\caption{The red squares represent the initial $\P$-positions for extension (b) of $(\frac{2}{7},\frac{10}{1})$-RN.}\label{F:5} 
\end{figure}

 The purpose of this paper has been to introduce the subject of $\B_Q$-subtraction games and resolve the most basic question, that is of $\varphi_Q$-equivalence of $Q$-subtraction games. As we indicate in Figures \ref{F:4} and \ref{F:5}, the most interesting ``new'' games are the generic $\B_Q$-subtraction games. (For example, could a study of extensions of $(\frac{1}{1},\frac{1}{0})\text{-\RW}$ lead to new revelations of the behavior of $(1,2)$GDWN?) Of course one would like to extend the notion of $\varphi_Q$-equivalence to games on several heaps, which can be done essentially by applying the same ideas as in Section \ref{S:1} and \ref{S:2}. We do not know of any literature on $\B_Q$-subtraction games even on two heaps.\\

\noindent{\bf Acknowledgment.}  I thank Ragnar Freij for an interesting discussion regarding generalizations to several heaps.\\\\
 
\begin{figure}[ht!]
\begin{center}
\includegraphics[width=0.5\textwidth]{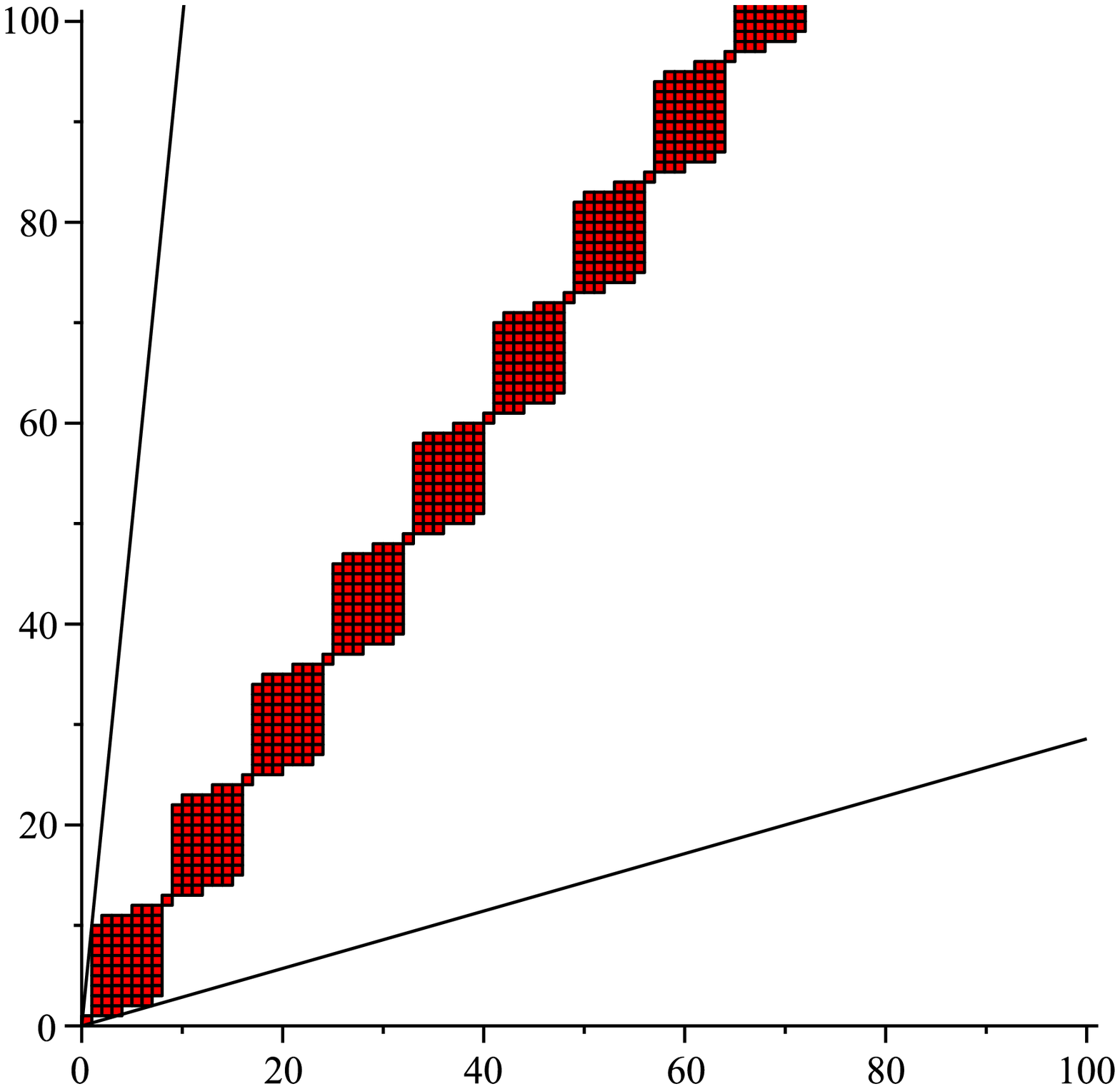}
\end{center}\caption{The red squares represent the initial $\P$-positions for extension (c) of $(\frac{2}{7},\frac{10}{1})$-RN.}\label{F:6} 
\end{figure}

\end{document}